\documentclass{article}
\usepackage{amsthm,amsmath,amssymb,amsthm,hyperref,xcolor}
\usepackage{graphicx,caption,appendix,multicol}
\usepackage[scale=0.8]{geometry}
\usepackage[numbers]{natbib}

\usepackage[utf8]{inputenc}

\newtheorem{prop}{Proposition}[]

\title{Asymptotic relatively more efficient test with auxiliary information: the case of the $Z$-test and the chi-square test}
\author{Mickael Albertus\footnote{mickael.albertus@gmail.com}}
\date{July 2019}

\begin{document}

\maketitle

\begin{multicols}{2}

\begin{abstract}
    The main goal of this article is to study how an auxiliary information can be used to improve the efficiency of two famous statistical tests: the $ Z$-test and the chi-square test. Many definitions of auxiliary information can be found in the statistical literature. In this article, the notion of auxiliary information is discussed from a very general point of view and depends on the relevant test. These two statistical tests are modified so that this information is taken into account. It is shown in particular that the efficiency of these new tests is improved in the sense of Pitman's ARE. Some statistical examples illustrate the use of this method.
\end{abstract}

\section{Introduction}

\subsubsection*{Main motivation} 
The main goal of this article is to present two new statistical tests which exploit a given auxiliary information. The new tests are based on modification of familiar statistical tests, the $Z$-test and the chi-square goodness-of-fit test and exploit a known auxiliary information in a way to get a more efficient test. These modifications are made so that, under the null hypothesis $ (H_0) $, the asymptotic behavior of the random variables involved by these test statistics does not change and, under $ (H_1) $, the probability of rejecting the null hypothesis is higher than that of the classical test. A description of the theoretical framework which allows for comparisons between two asymptotic statistical tests is provided below. To illustrate all results of this paper and to show how they can be used in a concrete way, these results will be applied with real data.

\subsubsection*{Auxiliary information} 
Although auxiliary information has been discussed extensively throughout the statistical literature, a generally-accepted definition does not exist. Statistical methods, such as stratification, calibration or Raking-Ratio, assume that a priori auxiliary information is given by the probability of sets of one or more partitions, and is hence known to the statistician. The following example illustrates this definition: a statistician working on a population of human beings knows the real proportion $p_M $ of men and $p_F $ of women in this population and wants to exploit this information in order to improve the estimates based on the sample taken from this population. If the statistician has a sample available, the frequency $ \widehat{p}_M $ of men and $ \widehat{p}_F $ of women in this sample is very likely to deviate from the known proportions $ p_M $ and $ p_F $. The methods cited above therefore aim to correct this difference with the hope of improving the estimates or the efficiency of statistical tests based on this sample. For instance, to estimate the unknown proportion $ p_S $ of sick people in the population, the raked estimator or Horvitz-Thompson estimator $$ \widehat{p}_S^\nabla = \widehat{p}_{M\cap S} \frac{p_M}{\widehat{p}_M} +\widehat{p}_{F\cap S} \frac{p_F}{\widehat{p}_F},  $$ can be considered, where $ \widehat{p}_{M \cap S} $ and $ \widehat{p}_{F  \cap S} $ are respectively the proportion of sick men and sick women in the population. This estimator is more precise than the frequency $ \widehat{p}_S $ of sick people.
Nevertheless, such a definition of the auxiliary information places constraints and is limited to only a handful of studies.
The current paper instead defines auxiliary information as information that asymptotically reduces variance of the main estimator implied by the test statistic. No additional assumptions are made about the nature and source of the information the statistician has at his disposal. A formal definition will be given later for each of the two tests and depends on that one.
Auxiliary information has been covered extensively in the statistical literature to study the estimators improved with additional information. To the author's knowledge, it has not been processed to study how to use the information to obtain more efficient tests.
The advantage of improving the efficiency of these tests is important since it makes it possible to accept smaller samples for a fixed level and power. The presented results are general and are applicable to many areas like medicine, biostatistics, economics and industry.

\subsubsection*{Asymptotic comparison of two tests} 
The asymptotic relative efficiency (ARE) plays the key role in this paper to compare different tests. More precisely, it is the Pitman relative efficiency which would be used to compare the new tests which exploit an auxiliary information versus the classical test which does not take into account this information. For two statistical tests ($i=1,2$) let $ n_i(\alpha, \beta,\theta) $ be the minimal sample sizes needed to test at a level $ \alpha $ and a power at least $ 1-\beta $ a null hypothesis $ (H_0):\theta=0 $ against a sequence of composite hypotheses $ (H_1):\theta =\theta_n $ with $ \theta_n $ a vanishing sequence. Notice that $ \theta $ is a parameter which can be a real or a real vector. The Pitman's ARE $ e_P $ of the first test with respect to the second test is defined, when this limit exists, as \begin{align}
 e_P = \lim_{\theta\to 0} \frac{n_2(\alpha,\beta,\theta)}{n_1(\alpha,\beta,\theta)}. \label{PitmanLimit}
\end{align} If this relative efficiency is smaller than 1, the second test needs a smaller size sample to attain the same level and power as the first test. In other words, the second test would be more efficient. As the same way, if this relative efficiency is larger than 1, the first test is then more effective. 

\subsubsection*{Framework} 
Let $ X_1,\dots,X_n, X $ be i.i.d. random variables defined on the same probability space $(\Omega,\mathcal{T},\mathbb{P})$ with same unknown distribution $ P=\mathbb{P}^{X} $ on some measurable space $ (\mathcal{X},\mathcal{T}') $. In order to get a probability space, the measurable space $(\mathcal{X},\mathcal{T}')$ is endowed with $P$. Let denote $ P(f) = \mathbb{E}[f(X)] $ and $ \mathbb{P}_n, \alpha_n $ respectively the empirical measure and process defined by \begin{align*}
    \mathbb{P}_n(f)&=\frac{1}{n} \sum_{i=1}^n f(X_i), \\
    \alpha_n(f) &= \sqrt{n}(\mathbb{P}_n(f)-P(f)),
\end{align*} for all $P$-measurable functions $f:\mathcal{X}\to \mathbb{R} $. By convenience, for $ A \in \mathcal{T}' $, $ \mathbb{P}_n(A)=\mathbb{P}_n(\mathbf{1}_A) $. For a $ P $-measurable function $f=(f_1,\dots,f_m):\mathcal{X}\to \mathbb{R}^m $ let denote \begin{align*}
     P[f] &= (P(f_1),\dots,P(f_m)),\notag \\
     \mathbb{P}_n[f] &= (\mathbb{P}_n(f_1),\dots,\mathbb{P}_n(f_m)).\notag 
\end{align*} The framework of this paper is non-parametric: no additional conditions are assumed on the law of $X$. The results are therefore applicable to a wide range of fields.

\subsubsection*{Organization} 
The new statistical tests which exploit the auxiliary information are presented and justified below. The two following sections describe the methods and main results for the new tests with auxiliary information. Section~\ref{Sec_ResultZTest} concerns the $Z$-test while Section~\ref{Sec_ResultChiTest} deals with the chi-square test. Section~\ref{Sec_Examp} provides examples for each of these improved tests as well as a non-exhaustive list of the literature surrounding the topic of auxiliary information.

\section{Pitman's ARE for the $Z$-test} \label{Sec_ResultZTest}

\subsubsection*{Notation} 
In this section, suppose that the random variables $ X_i $ are real. The common expectation and the variance of all variables $ X_i $ are respectively denoted by $ \mathbb{E}[X] $ and $ \sigma^2 $. This section focuses on the statistical $Z$-test based on the null hypothesis $ (H_0):\theta=0 $ with $ \theta=\mathbb{E}[X]-\mu $ and the alternative $ (H_1):\theta=h/\sqrt{n} $ for some $ h\in\mathbb{R} $. The classical statistic for this hypothesis is given by $$ Z_n = \frac{\sqrt{n}(\overline{X}_n-\mu)}{\widehat{\sigma}_n},$$ where $ \widehat{\sigma}_n $ is a consistent estimator of the standard deviation $ \sigma $, $\overline{X}_n=\frac{1}{n} \sum_{i=1}^n X_i $ is the empirical mean. Notice that if $ \sigma $ is a known value then $ \widehat{\sigma}_n=\sigma $ can be taken. Under $ (H_0) $, the statistic $ Z_n $ converges weakly to the normal distribution $ \mathcal{N}(0,1) $ while under $ (H_1) $ this statistic converges to $ \mathcal{N}\left( h/\sigma, 1 \right) $. Then the statistical test based on the rejection region $ |Z_n| > \Phi(1-\alpha/2) $ is an asymptotic confidence level $ \alpha $, where $ \Phi $ is the inverse of the standard normal cumulated distribution function.

\subsubsection*{Auxiliary information} 
In this context, an auxiliary information is an information which could be used to obtain an estimator of $ \mathbb{E}[X] $ with a lower variance than the natural empirical estimator $ \overline{X}_n $. To be in the most general framework, suppose that a known statistic $ \overline{X}_n^{\nabla} $ exploits an auxiliary information in the sense that it satisfies the following weak convergence \begin{align}
    \sqrt{n}(\overline{X}_n^{\nabla}-\mathbb{E}[X]) \underset{n \to +\infty}{\longrightarrow} \mathcal{N}(0,(\sigma^{\nabla})^2), \label{Cond0}
\end{align} where $ \sigma>\sigma^{\nabla} $. In other words, $ \overline{X}_n^\nabla $ is an estimator which uses the auxiliary information. Some examples of this known statistic and their associated value $ \sigma^{\nabla} $ for this test are given in Section~\ref{Sec_Examp}. The new statistic based on the $Z$-test with the auxiliary information is defined by $$ Z_n^{\nabla} = \frac{\sqrt{n}(\overline{X}_n^{\nabla}-\mu)}{\widehat{\sigma}^\nabla_n}, $$ where $ \widehat{\sigma}^\nabla_n $ is a consistent estimator of $ \sigma^{\nabla} $. As the same way as $ Z_n $, the new statistic $ Z_n^{\nabla} $ converges weakly to the normal distribution $ \mathcal{N}(0,1) $ under $ (H_0) $ and $ \mathcal {N}\left(h/\sigma^\nabla,1 \right) $ under $ (H_1) $. \medskip

\subsubsection*{Result} 
The following proposition suggests that the $ Z $-test is improved with the exploitation of an auxiliary information. This result is trivial since it is a direct application of the ARE definition.
\begin{prop}\label{PropZTest}
	The Pitman's ARE $ e_P $ of the classical test with respect to the new test which takes into account the auxiliary information satisfies $ e_P = (\sigma^\nabla/\sigma)^2 < 1 $.
\end{prop}
\begin{proof}
    It is an application of Theorem 14.19 of~\cite{Vaart1998AsymptoticStatistics} with $ \theta = \mathbb{E}[X]-\mu $, $ T_{n,1} = \overline{X}_n, M_1(\theta)=\theta, \sigma_1(\theta)=\widehat{\sigma}_n $ and $ T_{n,2} = \overline{X}_n^\nabla, M_2(\theta)=\theta, \sigma_2(\theta)=\widehat{\sigma}_n^\nabla $ which both satisfy Van der Vaart's condition (14.5).
\end{proof}

\noindent The interest of this proposition lies in its applications presented in Section~\ref{Sec_Examp}.

\section{Pitman's ARE for the goodness-of-fit test}\label{Sec_ResultChiTest}

\subsubsection*{Notation} 
In this section no assumption is made on the random variables $ X_i $, the distribution $ P $ or the set $ \mathcal{X} $. A parameter $ M \in \mathbb{N}\setminus \{0,1\} $ and a partition $ \mathcal{A} = (A_1,\dots, A_M) \subset \mathcal{T}' $ of $ \mathcal{X} $  such that $ P(A_i)\neq 0 $ for all $ i=1,\dots,M $ are fixed.  Let denote $ \mathcal{A}^*=(A_1,\dots,A_{M-1}) $ and remind that $ P[\mathcal{A}^*] $ and $ \mathbb{P}_n[\mathcal{A}^*] $ are the vectors respectively defined by \begin{align}
	P[\mathcal{A}^*] &= (P(A_1),\dots,P(A_{M-1}))\in \mathbb{R}^{M-1}, \notag \\
    \mathbb{P}_n[\mathcal{A}^*] &=  (\mathbb{P}_n(A_1),\dots,\mathbb{P}_n(A_{M-1})) \in \mathbb{R}^{M-1}. \label{DefPAStar}
\end{align}The goal of this test is to check if $ P[\mathcal{A}^*]=P_0[\mathcal{A}^*]=(P_0(A_1),\dots,P_0(A_{M-1})) $ for some measure $ P_0 $. The null hypothesis is $$ (H_0) : \Theta =\mathbf{0}_{M-1}, $$ where $ \mathbf{0}_{M-1}=(0,\dots,0)\in \mathbb{R}^{M-1} $ and \begin{align*}
	\Theta &= P[\mathcal{A}^*]-P_0[\mathcal{A}^*] \\
    &= \left(P(A_1)-P_0(A_1),\dots,P(A_{M-1})-P_0(A_{M-1}) \right).  
\end{align*} The simple sequence of alternative hypothesis considered for this test is $$ (H_1):\Theta=\mathbf{h}/\sqrt{n} $$ for some $ \mathbf{h} \in \mathbb{R}^{M-1} $. The chi-square test is based on the behavior of the random vector $\sqrt{n}(\mathbb{P}_n[\mathcal{A}^*]-P_0[\mathcal{A}^*]) $ which converges weakly under $(H_0)$ and $ (H_1) $ respectively to the multivariate normal law $ \mathcal{N}(0,\Sigma) $ and $ \mathcal{N}(\mathbf{h},\Sigma) $ where \begin{align}
	\Sigma =\mathrm{Diag}(P_0[\mathcal{A}^*]) - P_0[\mathcal{A}^*]^t \cdot P_0[\mathcal{A}^*]. \label{DefSig}
\end{align} According to Sherman–Morrison formula, $ \Sigma $ is invertible and $$ \Sigma^{-1} = \mathrm{Diag}\left(\frac{1}{P_0(A_1)},\dots,\frac{1}{P_0(A_{M-1})} \right)+\frac{1}{P_0(A_M)} \left( \begin{matrix}
    	1&\dots&1\\
        \vdots&\ddots&\vdots\\
        1&\dots&1
    \end{matrix}\right). $$  The statistic for the classic chi-square of goodness-of-fit test is given by \begin{align*} 
	\chi_n^2 &= n\sum_{i=1}^M \frac{(\mathbb{P}_n(A_i)-P_0(A_i))^2}{P_0(A_i)} \\
    &= n\sum_{i=1}^{M-1}  \frac{(\mathbb{P}_n(A_i)-P_0(A_i))^2}{P_0(A_i)} +\frac{n\left( \sum_{i=1}^{M-1} \mathbb{P}_n(A_i)-P_0(A_i)\right)^2}{P_0(A_M)} \\
    &=n(\mathbb{P}_n[\mathcal{A}^*]-P_0[\mathcal{A}^*]) \cdot \Sigma^{-1} \cdot (\mathbb{P}_n[\mathcal{A}^*]-P_0[\mathcal{A}^*])^t \\
    &= Z_n\cdot Z_n^t, \\
    Z_n &=\sqrt{n}(\mathbb{P}_n[\mathcal{A}^*]-P_0[\mathcal{A}^*])\cdot \Sigma^{-1/2}.
\end{align*}Under $ (H_0) $, $Z_n $ converges weakly to $ \mathcal{N}(\mathbf{0}_{M-1},\mathrm{Id}_{M-1}) $ and the statistic $ \chi_n^2 $ converges to $ \chi^2(M-1) $, a chi-square distribution with $M-1$ degrees of freedom. This proof is inspired by the first proof proposed in~\cite{Benhamou2018}. Other proofs of the convergence of the chi-squared statistics under $ (H_0) $ can be found in the last cited paper. Under $ (H_1) $, $Z_n$ converges to $ \mathcal{N}(\mathbf{M}, \mathrm{Id}_{M-1}) $ with $ \mathbf{M} = \mathbf{h}\cdot\Sigma^{-1/2} $ which leads to say that the statistic $ \chi_n^2 $ converges to a non-central chi-square distribution $ \chi^2(M-1; \lambda) $ -- see for example~\cite{Patnaik1949TheApplications,Cochran1952TheFit} -- with $ M-1 $ degrees of freedom and a non-centrality parameter \begin{align*}
	\lambda =\mathbf{M}\cdot \mathbf{M}^t = \mathbf{h}\cdot \Sigma^{-1}\cdot \mathbf{h}^t .
 \end{align*}

\subsubsection*{Auxiliary information} 
Suppose that an auxiliary information is available and that our aim is to taken into account this information to improve the chi-square test. Here, the auxiliary information is defined as an estimator $ \mathbb{P}_n^{(N)}[\mathcal{A}^*] $ of $ P[\mathcal{A}^*] $, given by~\eqref{DefPAStar}, with lower variance than the natural and empirical estimator $ \mathbb{P}_n[\mathcal{A}^*] $. Formally, the random vector $ \mathbb{P}_n^\nabla[\mathcal{A}^*] = \left( \mathbb{P}_n^\nabla(A_1),\dots,\mathbb{P}_n^\nabla(A_{M-1}) \right) \in \mathbb{R}^{M-1} $ is supposed to satisfy $$ \sqrt{n}(\mathbb{P}^{\nabla}_n[\mathcal{A}^*]-P[\mathcal{A}^*]) \underset{n\to+\infty}{\overset{\mathcal{L}}{\longrightarrow}} \mathcal{N}(0,\Sigma^{\nabla}), $$  where $ \Sigma^{\nabla} $ is a $ M\times M $ invertible covariance matrix such that \begin{align}
     \Sigma-\Sigma^{\nabla} \text{ is semi-definite positive}.\label{Cond1}
    \end{align} Condition~\eqref{Cond1} is what is called auxiliary information in this paper in the case of the chi-square test and this property will be essential for the next main result. Some examples of auxiliary information and matrices $ \Sigma^{\nabla} $ for this test which satisfy this hypothesis are given in Section~\ref{Sec_Examp}. Consistent estimator $  \widehat{\Sigma}^{\nabla}_{n} $ of $ \Sigma^{\nabla} $, like its empirical estimator, is considered. The chi-square statistic with auxiliary information is defined by \begin{align*}
    	\chi_n^{\nabla 2} &= n(\mathbb{P}_n[\mathcal{A}^*]-P_0[\mathcal{A}^*]) \cdot(\widehat{\Sigma}_n^\nabla)^{-1}\cdot (\mathbb{P}_n[\mathcal{A}^*]-P_0[\mathcal{A}^*])^t  \\
        &=Z_n^{\nabla} \cdot (Z_n^{\nabla})^t, \\
        Z_n &= \sqrt{n}(\mathbb{P}^{\nabla}_n[\mathcal{A}^*]-P_0[\mathcal{A}^*])\cdot (\widehat{\Sigma}_{n}^\nabla)^{-1/2}.
    \end{align*} Define the chi-square statistics with auxiliary information by matricially multiplying $ \sqrt{n}(\mathbb{P}_n^\nabla[\mathcal{A}^*]-P_0[\mathcal{A}^*]) $ by $ (\widehat{\Sigma}_n^\nabla)^{-1/2} $ for the definition of $Z_n^\nabla $ is motivated by the fact that $ Z_n^\nabla $ and therefore $ \chi_n^{2\nabla} $ follow the same law as $ Z_n $ and $ \chi_n^2 $ under $ (H_0) $. The statistical test based on the rejection region $ \chi_n^{\nabla 2} > t $, for some $t>0$, has the same alpha risk than the classical chi-square test based on the decision $ \chi_n^2 >  t $. Under $ (H_1) $, the random vector $ Z_n^\nabla $ converges weakly to $ \mathcal{N}(\mathbf{M}^\nabla,\mathrm{Id}_{M-1}) $ where $  \mathbf{M}^\nabla = \mathbf{h} \cdot (\Sigma^\nabla)^{-1/2} $ while $ \chi_n^{\nabla 2} $ converges to the non-central chi-square distribution $ \chi^2(M-1;\lambda^\nabla) $ with $ M-1$ degrees of freedom and the non-centrality parameter \begin{align*}
    	\lambda^\nabla = \mathbf{M}^\nabla\cdot (\mathbf{M}^\nabla)^t =  \mathbf{h} \cdot (\Sigma^\nabla)^{-1}\cdot\mathbf{h}^t.
\end{align*} If condition~\eqref{Cond1} is satisfied then $ \lambda \leqslant \lambda^\nabla $. The next paragraph shows that the efficiency of the chi-square test is increased when the auxiliary information is used.

\subsubsection*{Result} 
The main result concerning the efficiency of the chi-square test with auxiliary information is given by the following proposition.
\begin{prop}\label{PropChi2Stat}
The Pitman's ARE of the classical chi-square test with respect to the new chi-square test which takes into account the auxiliary information is given by
\begin{align*}
    e_P = \frac{\mathbf{h} \cdot \Sigma^{-1} \cdot \mathbf{h}^t}{\mathbf{h} \cdot (\Sigma^\nabla)^{-1} \cdot \mathbf{h}^t}.
\end{align*} This efficiency is bounded by $$ e_P \leqslant \lambda_{\max}((\Sigma^\nabla)^{-1} \Sigma) \leqslant 1, $$ where $\lambda_{\max}(\cdot) $ is the largest eigenvalue of a matrix.
\end{prop}
\noindent This proposition bounds to the Pitman's ARE of the new test and suggests that this one is more efficient than the classical chi-square test.
\begin{proof}
The minimal sample size needed to attain the level $ \alpha $ and the power $ 1-\beta $ for the classical chi-square test and the new chi-square test with auxiliary information are respectively denoted $ n_1(\alpha,\beta,\Theta) $ and $ n_2(\alpha,\beta,\Theta) $. The power of the tests without and with auxiliary information are respectively given by $ \pi_n,\pi_n^\nabla $ where \begin{align*}
        \pi_n &= \mathbb{P}\left(Z_n \cdot Z_n^t > \mathcal{Q}_{M-1}(\alpha) | H_1\right),\\
        \pi_n^\nabla &= \mathbb{P}\left(Z_n^\nabla \cdot (Z_n^\nabla)^t > \mathcal{Q}_{M-1}(\alpha) | H_1\right),
    \end{align*} where $\mathcal{Q}_M(t) $ is the $ t $-quantile of the $ \chi_M^2 $ distribution, that is \begin{align*}
    \mathcal{Q}_M(t)=\inf\{x : t \leqslant \mathbb{P}(X \leqslant x)\},
\end{align*} for a chi-square variable $  X \sim \chi^2(M-1) $. These powers can be approximated thanks to the non-central approximation:\begin{align*}
	\pi_n &= 1- F_{M-1}(\mathcal{Q}_{M-1}(\alpha)) + o(1)\\
    &= Q_{(M-1)/2}\left(\sqrt{\lambda_n}, \sqrt{\mathcal{Q}_{M-1}(\alpha)}\right) + o(1),\\
    \pi_n^\nabla &=  1- F_{M-1}^\nabla(\mathcal{Q}_{M-1}(\alpha)) + o(1) \\
    &=  Q_{(M-1)/2}\left(\sqrt{\lambda_n^\nabla}, \sqrt{\mathcal{Q}_{M-1}(\alpha)}\right) + o(1),
\end{align*} where $ F_{M-1}, F_{M-1}^\nabla $ are respectively the distribution functions of the non-central chi-square distribution $ \chi^2(M-1;\lambda_n), \chi^2(M-1;\lambda_n^\nabla) $ with $ \lambda_n = n \Theta\cdot \Sigma^{-1}\cdot \Theta^t, \lambda_n^\nabla=n \Theta \cdot (\Sigma^\nabla)^{-1}\cdot \Theta^t $, $ Q_{(M-1)/2} $ is the Marcus-$Q$-function and $ o(1) $ are sequences vanishing when $n\to+\infty $. Sequence of powers $ \pi_n,\pi_n^\nabla $ satisfy $ \pi_n\to1-\beta $ and $ \pi_n^\nabla\to 1-\beta $ if and only if $ \lambda_n \to G(1-\beta) $ and $ \lambda_n^\nabla \to G(1-\beta) $ when $ \Theta\to\mathbf{0}_{M-1}$ and $ G $ denoting the reciprocal of the application $$ x \mapsto Q_{(M-1)/2}\left(\sqrt{x}, \sqrt{\mathcal{Q}_{M-1}(\alpha)}\right). $$ This statement implies that $$ \lim_{\Theta\to\mathbf{0}_{M-1}} \frac{n_2(\alpha,\beta,\Theta) }{n_1(\alpha,\beta,\Theta) } \times \frac{\Theta \cdot (\Sigma^\nabla)^{-1} \cdot \Theta^t}{\Theta \cdot \Sigma^{-1} \cdot\Theta^t} = 1, $$ and consequently, \begin{align*}
	e_P&=\lim_{\Theta\to\mathbf{0}_{M-1}} \frac{n_2(\alpha,\beta,\Theta)}{n_1(\alpha,\beta,\Theta)} \\
	&= \lim_{\Theta\to\mathbf{0}_{M-1}} \frac{\Theta \cdot \Sigma^{-1} \cdot\Theta^t}{\Theta \cdot (\Sigma^\nabla)^{-1} \cdot \Theta^t} \\
	&=\frac{\mathbf{h} \cdot \Sigma^{-1} \cdot \mathbf{h}^t}{\mathbf{h} \cdot (\Sigma^\nabla)^{-1} \cdot \mathbf{h}^t}.
\end{align*} Since \begin{align*}
    \frac{\mathbf{h} \cdot \Sigma^{-1} \cdot \mathbf{h}^t}{\mathbf{h} \cdot (\Sigma^\nabla)^{-1} \cdot \mathbf{h}^t} = \frac{\mathbf{h}\cdot (\Sigma^\nabla)^{-1/2} \Sigma (\Sigma^\nabla)^{-1/2})\cdot \mathbf{h}^t}{\mathbf{h}\cdot \mathbf{h^t}},
\end{align*} then by Rayleigh-Ritz theorem, \begin{align*}
    e_P &\leqslant \max_{\mathbf{x}\neq \mathbf{0}} \frac{\mathbf{x}\cdot \left[(\Sigma^\nabla)^{-1/2} \Sigma (\Sigma^\nabla)^{-1/2}\right]\cdot \mathbf{x}^t}{\mathbf{x}\cdot \mathbf{x^t}} \\
    &\leqslant \lambda_{\max}((\Sigma^\nabla)^{-1/2} \Sigma (\Sigma^\nabla)^{-1/2})= \lambda_{\max}((\Sigma^\nabla)^{-1} \Sigma).
\end{align*} Condition~\eqref{Cond1} implies that $ \lambda_{\max}((\Sigma^\nabla)^{-1} \Sigma) \leqslant 1 $.
\end{proof}
\noindent Notice that if $ M=2 $ then $ e_P = \Sigma^{\nabla}/\Sigma $.

\section{Statistical examples}\label{Sec_Examp}

\subsubsection*{Organization} 
This section describes and justifies two methods, the Raking-Ratio method and the general auxiliary information, allowing to obtain an asymptotic reduction of variance and therefore an auxiliary information as defined in this paper. Other methods can also be used to obtain a variance reduction. For example, from a sufficient statistic $S $ the Rao-Blackwell theorem allows from an estimator $ \widehat{\theta} $ to construct a more precise estimator $ \mathbb{E}[\widehat{\theta}|S] $ of $ \theta $. This method is not detailed but suggests to use as auxiliary information $ \overline{X}_n^\nabla = \mathbb{E}[\overline{X}_n|S] $ for the $Z$-test and $ \mathbb{P}_n^\nabla[\mathcal{A}^*] = \mathbb{P}_n[\mathcal{A}^*|S]=(\mathbb{P}_n(A_1|S),\dots,\mathbb{P}_n(A_{M-1}|S))$ for the chi-square test with a sufficient statistic $ S $.  Subsection~\ref{SubSec_RR} deals with the Raking-Ratio method, that is a method which takes into account an auxiliary information given by the probabilities of sets of given partitions. Subsection~\ref{SubSec_InfauxGen} deals with the general auxiliary information. Examples and numerical simulations are given in these two subsections.

\subsubsection*{Framework} 
For all numerical examples presented in this section the distribution $ P = \mathbb{P}^X $ described by the Figure~\ref{LawXExample} is used.

\medskip
\noindent
    \includegraphics[scale=0.6]{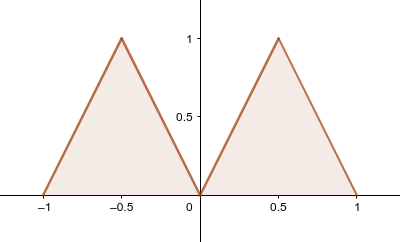}
    \captionof{figure}{Distribution of $ X $}
    \label{LawXExample}
\medskip
In particular, $ \mathbb{E}[X] = 0, \sigma^2=\mathrm{Var}(X)=7/24 \simeq 0.292 $.

\subsection{Raking-Ratio method}\label{SubSec_RR}

\subsubsection*{Presentation} 
The Raking-Ratio method is a statistical and computational method in order to take into account the auxiliary information given by the knowledge of $ P[\mathcal{A}^{(N)}] = (P(A_1^{(1)}), \dots,P(A_{m_N}^{(N)})) $ where $ \mathcal{A}^{(N)} = \{A_1^{(N)},\dots,A_{m_N}^{(N)}\} \subset \mathcal{T}' $ are partitions of $ \mathcal{X} $ and $ m_N \in\mathbb{N}^* $ is the size of the partition $ \mathcal{A}^{(N)} $, for all $ N \in\mathbb{N}^* $. It was introduced by Deming and Stephan~\cite{Deming1940} as a method to determine the projection of the empirical measure on the set of discrete probability measures satisfying all the constraints given by the auxiliary information. They did a mistake~\cite{Stephan1942} since this algorithm does not find this projection with respect to the chi-square distance. Ireland and Kullback~\cite{Ireland1968} have nevertheless proved that this method converges to this projection with respect to the Kullback-Leibler distance. This method was studied along the twentieth century -- see~\cite{G.J.BRACKSTONE1979,Sinkhorn1964,Konijn1981,Binder1988,Deville1993} -- but few of these results are satisfactory since they do not allow to prove theoretically that this method makes it possible to improve empirical estimates. Albertus and Berthet~\cite{Albertus2019b} studied this method from a new point of view, that of the empirical process theory, and they gave general results for infinite collection of estimators satisfying some metric entropy conditions. In~\cite{Albertus2019}, the author studied the behavior of the empirical process associated to the Raking-Ratio method when the information is given by an estimation or a learning and no more necessarily by an exact information and gave sufficient condition on this estimation to ensure to keep the same weak convergence. The two following paragraphs apply these theoretical  results in the case of the improved $Z$-test and chi-square test.

\subsubsection*{Raking-Ratio for the $Z$-test} 
\textbf{General case.} Result of Section~\ref{Sec_ResultZTest} can be applied if an empirical estimator of $ \mathbb{E}[X] $ satisfying condition~\eqref{Cond0} is known. The Raking-Ratio method gives a better estimator by exploiting iteratively the auxiliary information given by the knowledge of all $ P[\mathcal{A}^{(N)}] $. In our case, the raked empirical mean is given by $$ \overline{X}_n^{(N)} = \sum_{i=1}^{n} q_{n,i}^{(N)} X_i, $$ where $ q_{n,i}^{(N)} $ is the weight of $ X_i $ for the $ N $-th iteration of the Raking-Ratio method, that is, for all $ 1 \leqslant i \leqslant n $, $ q_{n,i}^{(0)} = 1/n $ and for $ N\in \mathbb{N} $, $$ q_{n,i}^{(N+1)} =q_{n,i}^{(N)} \left(\sum_{j=1}^M  \frac{P(A_j^{(N+1)})  \mathbf{1}_{A_j^{(N+1)}}(X_i) }{\sum_{k=1}^n q_{n,k}^{(N)} \mathbf{1}_{A_j^{(N+1)}}(X_k)}\right). $$The factor in bracket can be interpreted as corrections which operates the auxiliary information. An example of calculation of the raked empirical mean is given at the appendix A of~\cite{Albertus2019}. The asymptotic variance of $ \sqrt{n}\overline{X}_n^{(N)} $ is denoted by $ (\sigma^{(N)})^2  $. Albertus and Berthet proved that condition~\eqref{Cond0} is satisfied by taking $ X_n^{\nabla} = \overline{X}_n^{(N)} $ for some fixed $ N \in\mathbb{N} $, since they established that $$ (\sigma^{(N)})^2 = \sigma^2 - \sum_{k=1}^N (\Phi_k^{(N)})^t \cdot C_k \cdot \Phi_k^{(N)}, $$ where $ C_k  \in\mathcal{M}_{m_k,m_k} $ and $ \Phi_k^{(N)} \in \mathcal{M}_{m_k,1} $ are the matrix and the vector defined respectively by \begin{align}
    C_k &= \mathrm{Diag}(P[\mathcal{A}^{(k)}])-P[\mathcal{A}^{(k)}]^t \cdot P[\mathcal{A}^{(k)}], \label{DefCk} \\
    \Phi_k^{(N)} &=  \sum_{\substack{1\leqslant
L\leqslant N-k\\k<l_{1}<l_{2}<\dots<l_{L}\leqslant N}}(-1)^{L}\mathbf{P}%
_{\mathcal{A}^{(l_{1})}|\mathcal{A}^{(k)}}\mathbf{P}_{\mathcal{A}%
^{(l_{2})}|\mathcal{A}^{(l_{1})}}\notag\\
&\dots \mathbf{P}_{\mathcal{A}%
^{(l_{L})}|\mathcal{A}^{(l_{L-1})}} \left( \begin{smallmatrix}
    \mathbb{E}[X|A_1^{(l_L)}] \\
    \vdots \\
    \mathbb{E}[X|A_{m_{l_L}}^{(l_L)}]
\end{smallmatrix} \right) + \left( \begin{smallmatrix}
    \mathbb{E}[X|A_1^{(k)}] \\
    \vdots \\
    \mathbb{E}[X|A_{m_{k}}^{(k)}]
\end{smallmatrix} \right), \notag
\end{align} and $ \mathbf{P}_{\mathcal{A}^{(i)}|\mathcal{A}^{(j)}} \in\mathcal{M}_{m_j,m_i} $ are stochastic matrices defined for all $ i,j\in \mathbb{N}^*  $ by \begin{align}
     (\mathbf{P}_{\mathcal{A}^{(i)}|\mathcal{A}^{(j)}})_{k,l} = P(A_l^{(i)}|A_k^{(j)}), \label{DefMatStoRR}
\end{align} for all $ 1 \leqslant l \leqslant m_i $ and $ 1 \leqslant k \leqslant m_j $. Since $ C_k $ are covariance matrices, and in particular semi-definite positive matrices, then $ \sigma^{(N)} \leqslant \sigma $ for all $ N \in \mathbb{N} $. Notice that these last matrices depend only on the auxiliary information given by all $ P[\mathcal{A}^{(N)}] $.\medskip

\textbf{Simple case.} In~\cite{Albertus2019} the author gave, in a general way, some examples of possible and explicit values of $ \sigma^{(N)} $ for $ N=1,2 $ when the simple case $ \mathcal{A}^{(1)}=\{A,A^C\},\mathcal{A}^{(2)}=\{B,B^C\} $ is considered. In our case these values are \begin{align*}
    (\sigma^{(1)})^2 &= \sigma^2 - \frac{p_A}{p_{\overline{A}}} \Delta_A^2,\\
    (\sigma^{(2)})^2&=\sigma^2-\frac{p_B}{p_{\overline{B}}} \Delta_B^2  -K \Delta_A^2, \\
\end{align*} where $ \Delta_A = \mathbb{E}[X|A]-\mathbb{E}[X] $ and $ \Delta_B = \mathbb{E}[X|B]-\mathbb{E}[X] $ and \begin{align*}
    p_A&=P(A),\quad p_{\overline{A}} = P(A^C), \\
    p_B&=P(B),\quad p_{\overline{B}} = P(B^C), \quad p_{A \cap B} = P(A\cap B),\\
    K &= p_A p_{\overline{A}}+\frac{p_B p_{\overline{B}}(p_{A\cap B}-p_A p_B)}{p_A^2 p_{\overline{A}}^2}.
\end{align*} If $ \mathcal{A}^{(2k-1)} = \mathcal{A}^{(1)} $ and $ \mathcal{A}^{(2k)} = \mathcal{A}^{(2)} $ for $ k > 1 $ then the value $ \sigma^{(\infty)} $, the standard deviation of $ \overline{X}_n^{(N)} $ when $ N $ goes to infinity, that is  when the Raking-Ratio method converges, is given by the following formula -- see (2.12) of~\cite{Albertus2019}: \begin{align*} 
	(\sigma^{(\infty)})^2 &= \sigma^2 - \frac{K p_A p_B}{p_A p_B p_{\overline{A}} p_{\overline{B}}-(p_{AB}-p_A p_B)^2},  \\
    K&=  p_A \Delta_A^2+p_B \Delta_B^2-p_A p_B (\Delta_A-\Delta_B)^2\\
    &\qquad-2p_{AB} \Delta_A \Delta_B.
\end{align*}  When the events $ A $ and $ B $ are independent then \begin{align*}
    (\sigma^{(\infty)})^2 = (\sigma^{(2)})^2 = (\sigma^{(1)})^2-\frac{p_B}{p_{\overline{B}}} \Delta_B^2.
\end{align*} In the independent case, since $ \sigma^{(\infty)} = \sigma^{(2)} $, the Raking-Ratio method could be stopped with $ N = 2 $ steps.\medskip

\textbf{Numerical simulation.} The previous results are applied with $ X $ following the distribution given by Figure~\ref{LawXExample} and the following  independent sets \begin{align}
    A = \{X \in [-0.5,0]\cup[0.5,1]\}, \quad B &= \{  X \leqslant 0 \} \label{DefABExamples}
\end{align} which satisfy \begin{align*}
    p_A=p_{\overline{A}} = p_B=p_{\overline{B}}&=0.5, \\
    \mathbb{E}[X|A] &= 1/6, \\
    \mathbb{E}[X|B]&=-0.5. 
\end{align*} With these sets, the empirical estimator with auxiliary information is given by $ \overline{X}_n^{(N)} = \sum_{i=1}^n q_{n,i}^{(N)} X_i $ where $ q_{n,i}^{(0)}=1/n $ and for $ N \in \mathbb{N} $, \begin{align*} 
	q_{n,i}^{(2N+1)} &= \frac{q_{n,i}^{(2N)} }{2} \times\frac{\mathbf{1}_A(X_i)}{\sum_{k=1}^n q_{n,k}^{(2N)} \mathbf{1}_A(X_k)} \\
    &\quad + \frac{q_{n,i}^{(2N)} }{2}\times \frac{\mathbf{1}_{\overline{A}}(X_i)}{\sum_{k=1}^n q_{n,k}^{(2N)} \mathbf{1}_{\overline{A}}(X_k)}, \\
	q_{n,i}^{(2(N+1))} &= \frac{q_{n,i}^{(2N+1)} }{2} \times \frac{\mathbf{1}_B(X_i)}{\sum_{k=1}^n q_{n,k}^{(2N+1)} \mathbf{1}_B(X_k)}\\
    &\quad + \frac{q_{n,i}^{(2N+1)} }{2}\times  \frac{\mathbf{1}_{\overline{B}}(X_i)}{\sum_{k=1}^n q_{n,k}^{(2N+1)} \mathbf{1}_{\overline{B}}(X_k)}. 
\end{align*}For example for $ N=1 $,\begin{align*}
	\overline{X}_n^{(1)} &= \frac{1}{2}\left(\frac{\sum_{i=1}^n X_i \mathbf{1}_{X_i \in A}}{\sum_{i=1}^n \mathbf{1}_{X_i \in A}}+\frac{\sum_{i=1}^n X_i \mathbf{1}_{X_i \in \overline{A}}}{\sum_{i=1}^n \mathbf{1}_{X_i \in \overline{A}}} \right).
\end{align*} The asymptotic variances of $ \overline{X}_n^{(N)} $ for $ N=1 $ and $N=2 $ or $N=\infty $ are equal to \begin{align*} 
    (\sigma^{(1)})^2 &= \sigma^2-1/36 = 19/72 \simeq 0.264, \\
    (\sigma^{(\infty)})^2 &=(\sigma^{(2)})^2 = \sigma^2-1/4 = 1/24 \simeq 0.042.
\end{align*} Figure~\ref{distribution_ztest_RR_H0} represents the distribution of $ \sqrt{n} \overline{X}_n $ and $ \sqrt{n}\overline{X}_n^\nabla = \sqrt{n} \overline{X}_n^{(N)} $ for $ N=1, 2 $ which are close to $ \mathcal{N}(0,\sigma^2) $ and $ \mathcal{N}(0,(\sigma^{(N)})^2)$. The decrease in variance is particularly visible for $ N=2 $ or $ N=\infty $. 

\includegraphics[width=250pt]{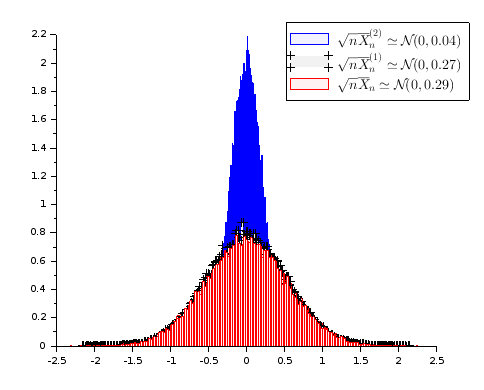}
  \captionof{figure}{Distribution of $ \sqrt{n}\overline{X}_n $ and $ \sqrt{n}\overline{X}_n^{\nabla} $ for $ N=1,2 $}\label{distribution_ztest_RR_H0}
  According to Proposition~\ref{PropZTest}, the Pitman's ARE for $ N = 1 $ is equal to $$ e_P^{(1)} = (\sigma^{(1)}/\sigma)^2=19/21 \simeq 0.905.  $$ For $ N = 2 $ or $ N=\infty $, \begin{align*}
    e_P^{(2)}&= (\sigma^{(2)}/\sigma)^2 = 1/7\simeq 0.143 \\
    e_P^{(\infty)} &= (\sigma^{(\infty)}/\sigma)^2=e_P^{(2)}.
\end{align*} Figure~\ref{distribution_ztest_RR_H1} represents the distribution of $ Z_n $ and $ Z_n^{\nabla} $ when $ \theta=0.5/\sqrt{n} $ and $n=100$. This figure illustrates that $ Z_n, Z_n^{(1)}, Z_n^{(2)} $ are asymptotically close to respectively $ \mathcal{N}(-\sqrt{6/7},1),  \mathcal{N}(-\sqrt{18/19},1) $ and $  \mathcal{N}(-\sqrt{6},1) $. The case $ N=2 $ or $ N=\infty $ is the most interesting because what makes the tests with auxiliary information more effective is highlighted, that is to say that the expectation of $Z_n^\nabla $ takes expected values greater than that of $ Z_n $.
    \includegraphics[width=260pt]{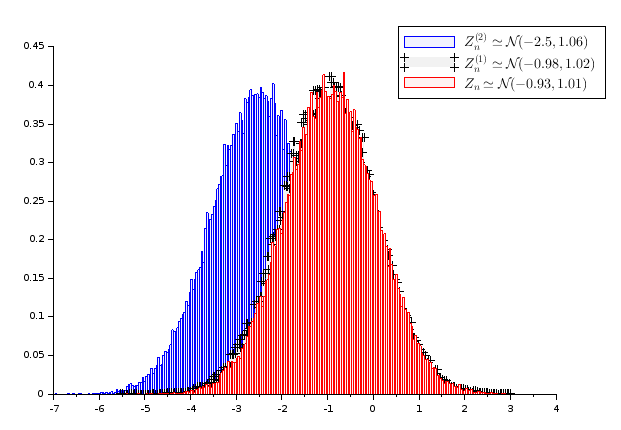}
  \captionof{figure}{Distribution of $ Z_n $ and $ Z_n^{\nabla} $ for $ N=1,2 $ and $n=100$ under $ (H_1) $}\label{distribution_ztest_RR_H1}

\subsubsection*{Raking-Ratio for the chi-square tests}
\textbf{General case.} To use the result of Section~\ref{Sec_ResultChiTest}, an estimator of $ P[\mathcal{A}^*] $ more efficient than the empirical estimator $ \mathbb{P}_n[\mathcal{A}^*]  $ in the sense given by condition~\eqref{Cond1} must be known to the statistician. The Raking-Ratio algorithm gives again a better estimator with the auxiliary information by using the knowledge of all $ P[\mathcal{A}^{(N)}] $ for $ N \geqslant 1 $. The raked estimator of $ P[\mathcal{A}^*] $ is  $$ \mathbb{P}_n^{(N)}[\mathcal{A}^*] = \left(\mathbb{P}_n^{(N)}(A_1), \dots,\mathbb{P}_n^{(N)}(A_{M-1})\right), $$ where $ \mathbb{P}_n^{(N)}(A) $ is iteratively defined, for a measurable set $ A $, by $ \mathbb{P}_n^{(0)}(A) = \mathbb{P}_n(A) $ and for all $ N \in \mathbb{N} $, $$ \mathbb{P}_n^{(N+1)}(A) = \sum_{j=1}^{m_{N+1}} \frac{\mathbb{P}_n^{(N)}(A \cap A_j^{(N+1)}) P(A_j^{(N+1)})}{\mathbb{P}_n^{(N)}(A_j^{(N+1)})}. $$ Results of Albertus and Berthet imply in particular that the process $$ \alpha_n^{(N)}[\mathcal{A}^*] = \sqrt{n}\left(\mathbb{P}_n^{(N)}[\mathcal{A}^*]-P[\mathcal{A}^*]\right), $$ converges to the singular multivariate normal distribution $ \mathcal{N}(0,\Sigma^{(N)})$ with $ \Sigma^{(N)} $ is the covariance matrix defined by \begin{align} 
	\Sigma^{(N)}=\Sigma - \sum_{k=1}^N (\Phi_{k}^{(N)})^t \cdot C_k \cdot \Phi_{k}^{(N)}  \label{DefSigNRRChiSquare}
\end{align} where $ \Sigma $ is defined by~\eqref{DefSig}, $ C_k  \in\mathcal{M}_{m_k,m_k} $ are the same covariance matrices defined above by~\eqref{DefCk} and vectors $ \Phi_k^{(N)} \in \mathcal{M}_{m_k,M-1} $ are the matrices whose the $i$th column is given for all by \begin{align*}
    &(\Phi_{k}^{(N)})_{\cdot,i} =\left( \begin{smallmatrix}
    P(A_i|A_1^{(k)}) \\
    \vdots \\
    P(A_i|A_{m_{k}}^{(k)})
\end{smallmatrix} \right)\\
    &+\sum_{\substack{1\leqslant
L\leqslant N-k\\k<l_{1}<l_{2}<\dots<l_{L}\leqslant N}}(-1)^{L}\mathbf{P}%
_{\mathcal{A}^{(l_{1})}|\mathcal{A}^{(k)}}\mathbf{P}_{\mathcal{A}%
^{(l_{2})}|\mathcal{A}^{(l_{1})}}\\
&\dots 
\mathbf{P}_{\mathcal{A}%
^{(l_{L})}|\mathcal{A}^{(l_{L-1})}} \left( \begin{smallmatrix}
    P(A_i|A_1^{(l_L)}) \\
    \vdots \\
    P(A_i|A_{m_{l_L}}^{(l_L)})
\end{smallmatrix} \right), \\
    &(\Phi_{N}^{(N)})_{\cdot,i} = \left( \begin{smallmatrix}
    P(A_i|A_1^{(N)}) \\
    \vdots \\
    P(A_i|A_{m_{N}}^{(N)})
\end{smallmatrix} \right),
\end{align*} for all $ i=1,\dots,M-1 $, $ N \geqslant 1 $ and $ 1 \leqslant k < N $ and $ \mathbf{P}_{\mathcal{A}^{(i)}|\mathcal{A}^{(j)}} \in\mathcal{M}_{m_j,m_i} $ are stochastic matrices defined by~\eqref{DefMatStoRR}. The matrix $ \Sigma^{(N)} $ ensures condition~\eqref{Cond1} so the results of Proposition~\ref{PropChi2Stat} can be applied with $ \mathbb{P}_n^\nabla[\mathcal{A}^*]=\mathbb{P}_n^{(N)}[\mathcal{A}^*] $ for any $ N \geqslant 1 $.\medskip

\textbf{Simple case.} The simple case with the following parameters is studied. The auxiliary information is given by the knowledge of $ P[\mathcal{A}^{(1)}] $ and $ P[\mathcal{A}^{(2)}] $ with $ \mathcal{A}^{(1)} = \{A_1, A_1^C\}, \mathcal{A}^{(2)}=\{A_2, A_2^C\} $. The categories for the chi-square test are given by $ M=2 $ and $ \mathcal{A} = \{A, A^C\} $. Matrices $ \Sigma^{(1)} $ and $ \Sigma^{(2)} $ require calculating $ \Phi_1^{(1)}, \Phi_1^{(2)} $ and $ \Phi_2^{(2)} $ which are in this case equal to \begin{align*}
    &\Phi_1^{(1)} =  \left( \begin{matrix}
        P(A|A_1) \\ P(A|A_1^C)
    \end{matrix}\right) \\
    &\Phi_2^{(2)} =  \left( \begin{matrix}
        P(A|A_2) \\ P(A|A_2^C)
    \end{matrix}\right) \\
    &\Phi_1^{(2)} =\left( \begin{matrix}
        P(A|A_1) \\ P(A|A_1^C)
    \end{matrix}\right)  - \left( \begin{matrix}
       P(A_2|A_1) & P(A_2^C|A_1) \\
       P(A_2|A_1^C) & P(A_2^C|A_1^C)
    \end{matrix}\right) \cdot  \left( \begin{matrix}
        P(A|A_2) \\ P(A|A_2^C)
    \end{matrix}\right)
\end{align*}

\textbf{Numerical simulation.} Results are applied with $ \mathcal{A}=\{A,A^C\},  \mathcal{A}^{(1)}=\{A_1,A_1^C\}, \mathcal{A}^{(2)}=\{A_2,A_2^C\}  $ with \begin{align*}
     A &= \{ X \leqslant 0.5 \}, \\
     A_1&=\{X\in [-0.5,0]\cup [0.5,1]\}, \\
     A_2 &= \{X\leqslant 0\}.
\end{align*}  In particular $ P[\mathcal{A}^*]=P(A) = 3/4 $ and event $A$ is dependent from the events $ A^{(1)}, A^{(2)} $. With these values, matrices $ \Sigma $ and $ \Sigma^{(N)} $ are real values and  \begin{align*}
    \Sigma &= P(A)(1-P(A)) =3/16\\
    \Phi_1^{(1)} &= \left(\begin{matrix}
        1/2 \\ 1
    \end{matrix}\right), \\
    \Phi_1^{(2)} &= \left(\begin{matrix}
       -1/4 \\ 1/4
    \end{matrix}\right),\\
    \Phi_2^{(2)} &= \left(\begin{matrix}
       1\\ 1/2
    \end{matrix}\right), \\
    C_1 &= C_2 = \frac{1}{4}\left(\begin{matrix}
        1 & -1 \\
        -1 & 1
    \end{matrix} \right).
\end{align*} According to~\eqref{DefSigNRRChiSquare},
\begin{align*}
    \Sigma^{(1)} &= \Sigma - (\Phi_1^{(1)})^t \cdot C_1 \cdot \Phi_1^{(1)} = 1/8, \\
    \Sigma^{(2)} &= \Sigma - (\Phi_1^{(2)})^t \cdot C_1 \cdot \Phi_1^{(2)} - (\Phi_2^{(2)})^t \cdot C_2 \cdot \Phi_2^{(2)}=1/16.
\end{align*} Condition~\eqref{Cond1} is met as expected for $N=1,2 $. According to Proposition~\ref{PropChi2Stat}, with $ \mathbb{P}^{\nabla}_n[\mathcal{A}^*] = \mathbb{P}_n^{(N)}[\mathcal{A}^*]=\mathbb{P}_n^{(N)}(A) $, the Pitman's ARE for $ N=1 $ is equal to \begin{align*}
	e_P^{(1)} & = \Sigma^{(1)}/\Sigma=2/3,
\end{align*} and for  $ N=2 $,  $$ e_P^{(2)} = \Sigma^{(2)}/\Sigma= 1/3. $$ Figure~\ref{distribution_chi2test_RR_H0} illustrates the asymptotic distribution of $ \sqrt{n}(\mathbb{P}_n(A)-P(A)) $ and $ \sqrt{n}(\mathbb{P}_n^{(N)}(A)-P(A)) $ for $ N=1,2 $. 

\includegraphics[width=250pt]{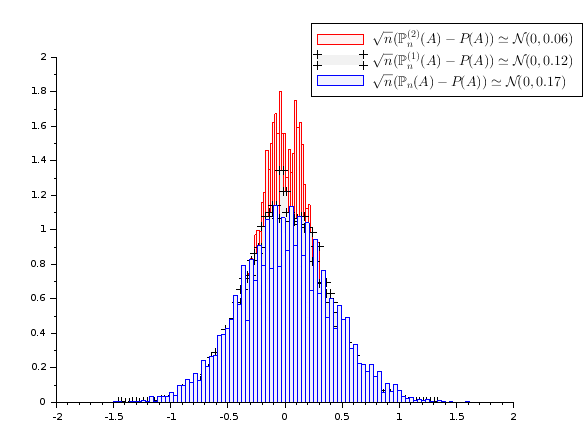}
  \captionof{figure}{Distribution of $ \sqrt{n}(\mathbb{P}_n(A)-P(A)) $ and $ \sqrt{n}(\mathbb{P}_n^{(N)}(A)-P(A)) $ for $ N=1,2 $}\label{distribution_chi2test_RR_H0}
  
The hypothesis $ \Theta =P(A)+0.5/\sqrt{n}= 3/4+0.5/\sqrt{n} $ is considered to be in the case of $ (H_1)$ hypothesis. Figure~\ref{distribution_chi2test_RR_H1} represents the distribution of $ \chi_n^2 $ and $ \chi_n^{\nabla 2} $ for $ N=1,2 $ which are respectively close to $ \chi^2(1;4/3), \chi^2(1;2) $ and $ \chi^2(1;4) $.

\includegraphics[width=225pt]{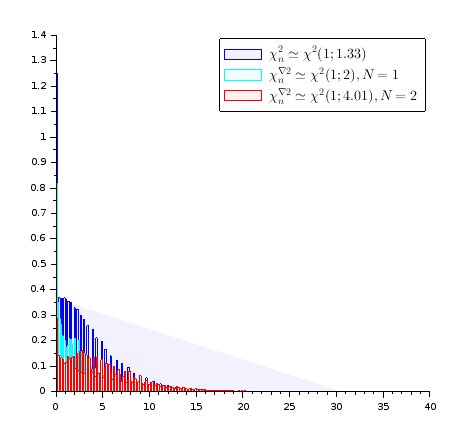}
  \captionof{figure}{Distribution of $ \chi_n^2 $ and $ \chi_n^{\nabla 2} $ for $ N=1,2 $ under $ (H_1) $ with $ n=100 $}\label{distribution_chi2test_RR_H1}

\subsection{General auxiliary information} \label{SubSec_InfauxGen}

\subsubsection*{Literature} 
Many definitions of the auxiliary information could be found among the literature concerning this subject. Some authors gave their own definition of the auxiliary information and they established some results in order to prove the efficiency of the use of auxiliary information. Nevertheless, most part of these definitions does not include some natural examples of what auxiliary information can be. For example, Dmitriev and Tarasenko~\cite{Dmitriev1992} defined the auxiliary information as the knowledge of probabilities $ P(g_1),\dots,P(g_m) $ for some measurable functions $ g_1,\dots,g_m $ or the knowledge of  an approximation of these probabilities. They determined the projection of the empirical measure which minimize the Kullback-Leibler divergence over the set of the probability measures satisfying the information. For another example, in~\cite{Zhang1995,Zhang1997,Zhang1997a}, Zhang defined the auxiliary information as a known function which the expectation cancels, that is a measurable real-valued function $ g $ such that $ \mathbb{E}[g(X)]=0 $ and he established some results when an information of this kind is known. However, these definitions do not scope the concept of general auxiliary information presented in the next paragraph. For instance, the knowledge of the variance of $ X $ can be not supported. As a matter of fact, $ \mathrm{Var}(X) = \mathbb{E}[h(X)] $ with $ h(X)=(X-\mathbb{E}[X])^2 $ implies that $ \mathbb{E}[g(X)] = 0 $ for $ g(X)=h(X)-\mathrm{Var}(X) $ and $ g $ can be not computable if $ \mathbb{E}[X] $ is unknown.

\subsubsection*{General auxiliary information} 
A very general definition is given by Tarima and Pavlov -- see~\cite{Tarima2006} -- where the auxiliary information is defined as an unbiased estimator which satisfies a CLT. Their result is general since the auxiliary information can be given by several sources of auxiliary information and can even be uncertain, that is an estimate of the true auxiliary information under certain conditions of asymptotic normality. To illustrate the results of Sections~\ref{Sec_ResultZTest} and~\ref{Sec_ResultChiTest}, the two following paragraphs present the case when the conditional mean of interest variables is known. More formally, suppose that $ X $ is real random variable such that its conditional expectation conditional on an event is known. In this paragraph, let work on real random variables $ X_1,\dots,X_n $ such that the expectation of $ X $ conditional on the event that $ X $ belongs to some predefined set $ C \in \mathcal{T}' $ is known. Thereby the auxiliary information is given by $ \mathbb{E}[X|C] = \mathbb{E}[X\mathbf{1}_C]/P(C) $. This kind of auxiliary information is not supported by the definitions of auxiliary information given by Dmitriev and Tarasenko or Zhang recalled below.

\subsubsection*{Conditional mean auxiliary information for the $Z$-test} 
\textbf{General case.} In this paragraph, the way the auxiliary information $ \mathbb{E}[X|C]$  can be exploited in the case of the $Z$-test is presented. The natural empirical estimator of the auxiliary information $ \mathbb{E}[X|C] $ is denoted by \begin{align}
    \mathbb{P}_n(X|C) = \frac{\mathbb{P}_n(X\mathbf{1}_C)}{\mathbb{P}_n(C)}. \label{DefPnX|C}
\end{align} The aim of this paragraph is to take into account this auxiliary information to suggest an estimator $ X_n^{\nabla} $ of $ \mathbb{E}[X] $ with a lower variance than the natural empirical estimator $ \overline{X}_n $. With this new estimator, condition~\eqref{Cond0} would be satisfied. To make the parallel with the article of Tarima and Pavlov, the same notation than their article is adopted. Suppose in this paragraph that one have one exact auxiliary information but an uncertain auxiliary information can be considered, given for example by an estimate based on another larger independent sample.  In our case there is $I=1 $ data source, $ J_1=1 $ auxiliary information and \begin{align*}
    \Theta &= \mathbb{E}[X], \quad \widehat{\Theta} = \overline{X}_n = \frac{1}{n} \sum_{i=1}^n X_i, \\
    \widetilde{\mathcal{B}} &= \mathcal{B} = \mathbb{E}[X|C], \quad \widehat{\mathcal{B}} = \mathbb{P}_n(X|C).
\end{align*} With these values, \begin{align*}
    K_{11} &= \mathrm{Var}(\widehat{\Theta}) = \sigma^2 \\
    K_{12} &= \mathrm{Cov}(\widehat{\mathcal{B}}, \widehat{\Theta}) = \mathrm{Cov}(\mathbb{P}_n(X|C), \overline{X}_n) \\
    K_{22}' &= \mathrm{Var}(\widehat{\mathcal{B}}) = \mathrm{Var}(\mathbb{P}_n(X|C)), \\ K_{22}''&=\mathrm{Var}(\widetilde{\mathcal{B}})=0, \\
    K_{22} &= K_{22}'+K_{22}'' = K_{22}' = \mathrm{Var}(\mathbb{P}_n(X|C)).
\end{align*} The elements $ K_{11}, K_{12}, K_{22} $ are unknown or could not be expressed simply since the estimator $ \mathbb{P}_n(X|C) $ is a quotient of the empirical measure. Therefore the first suggested estimator of Tarima and Pavlov \begin{align*}
     \widehat{\Theta}^0 &= \widehat{\Theta}-K_{12}K_{22}^{-1} (\widehat{\mathcal{B}}-\widetilde{\mathcal{B}}) \\
     &= \overline{X}_n-K_{12}K_{22}^{-1} (\mathbb{P}_n(X|C)-\mathbb{E}[X|C]),
\end{align*} which exploits the auxiliary information $ \mathbb{E}[X|C] $ is uncomputable. This case is common and it is for that reason that these authors suggested to replace these unknown values by consistent estimators of them. Values $ K_{11},K_{12},K_{22} $ can be estimated respectively by the following values \begin{align*}
    \widehat{K}_{11} &= \widehat{\sigma}_n^2,\\
    \widehat{K}_{12} &= \frac{1}{n} \left(\mathbb{P}_n(X^2|C)-\overline{X}_n \mathbb{E}[X|C] \right) , \\
    \widehat{K}_{22} &= \frac{1}{n\mathbb{P}_n(C)} \left(\mathbb{P}_n(X^2|C)-\mathbb{E}[X|C]^2 \right),
\end{align*} where $ \widehat{\sigma}_n^2 $ is a consistent estimator of $ \sigma^2 $.  With these consistent estimators, Tarima and Pavlov suggest to use the statistic \begin{align*}
    \widehat{\Theta}^* &= \widehat{\Theta} - \widehat{K}_{12} \widehat{K}_{22}^{-1} (\widehat{\mathcal{B}}-\widetilde{\mathcal{B}}) \\
    &=\overline{X}_n - \widehat{K}_{12} \widehat{K}_{22}^{-1} (\mathbb{P}_n(X|C)-\mathbb{E}[X|C]).
\end{align*} By taking $ a_n= \sqrt{n} $ the first Tarima and Pavlov conditions mentioned in Section 1.3 of their paper are respected. More precisely, $ \zeta_{1n}=0, \Sigma_{22}''=0 $ and \begin{align*}
    \xi_n = a_n(\widehat{\Theta}-\Theta)&\underset{n\to+\infty}{\overset{\mathcal{L}}{\longrightarrow}} \xi = \mathcal{N}(0,\Sigma_{11}), \\
    \tau_n = a_n (\widehat{\mathcal{B}}-\mathcal{B}) &\underset{n\to+\infty}{\overset{\mathcal{L}}{\longrightarrow}} \tau = \mathcal{N}(0,\Sigma_{22}'), \\
    n K_{22}' = \mathrm{Var}(\alpha_n(X|C))  &\underset{n \to +\infty}{\longrightarrow} \Sigma_{22}',
    \end{align*} with \begin{align*}
        \alpha_n(X|C)&=\sqrt{n}(\mathbb{P}_n(X|C)-\mathbb{E}[X|C]), \\ 
        \Sigma_{11}&=K_{11}=\sigma^2, \\
        \Sigma_{22}'&= \frac{\mathrm{Var}(X|C)}{P(C)},
    \end{align*} where $ \mathrm{Var}(X|C) $ is the variance of $ X $ conditional on the event $ C $ defined by  \begin{align}
    	\mathrm{Var}(X|C) &= \mathbb{E}[X^2|C]-\mathbb{E}[X|C]^2. \label{DefVar|C}
    \end{align} According to Proposition 1 of Tarima and Pavlov, \begin{align*}
        a_n(\widehat{\Theta}^0-\Theta)\underset{n \to +\infty}{\overset{\mathcal{L}}{\longrightarrow}} &\mathcal{N}\left(0,(\sigma^\nabla)^2\right), 
    \end{align*} where \begin{align*}
    (\sigma^\nabla)^2 &= \Sigma_{11}- \frac{\Sigma_{12}^2}{\Sigma_{22}}=\sigma^2-\frac{\mathrm{Cov}^2(X\mathbf{1}_C, X)}{P(C)\mathrm{Var}(X|C)}, \\
        \Sigma_{12} &= \mathrm{Cov}(\xi,\tau)= \frac{\mathrm{Cov}(X\mathbf{1}_C,X)}{P(C)},\\
        \Sigma_{22} &= \Sigma_{22}' = \frac{\mathrm{Var}(X|C)}{P(C)}. 
    \end{align*} Conditions of Proposition 2 of Tarima and Pavlov are respected since the following asymptotic behaviour is satisfied: \begin{align*}
    a_n^2(\widehat{K}_{12}-K_{12}) &\underset{n\to+\infty}{\overset{\textrm{a.s.}}{\longrightarrow}} 0, \\
    a_n^2(\widehat{K}_{22}-K_{22}) &\underset{n\to+\infty}{\overset{\textrm{a.s.}}{\longrightarrow}} 0.
\end{align*} By taking $ X_n^{\nabla} =\widehat{\Theta}^*  $, Proposition 2 of Tarima and Pavlov imply that \begin{align*}
    \sqrt{n}( X_n^{\nabla}-\mathbb{E}[X])\underset{n \to +\infty}{\overset{\mathcal{L}}{\longrightarrow}} &\mathcal{N}\left(0,(\sigma^\nabla)^2\right).
\end{align*} According to Proposition~\ref{PropZTest} of this article, the Pitman's ARE $ e_P $ is equal to $$ e_P = (\sigma^\nabla/\sigma)^2 = 1-\frac{\mathrm{Cov}^2(X\mathbf{1}_C,X)}{\sigma^2 P(C)\mathrm{Var}(X|C)}. $$  \textbf{Numerical simulation.} If the previous results are applied with $ C = \{|X| \leqslant 0.5\} $, that is the auxiliary information is given by the knowledge of the value $ \mathbb{E}[X| \ |X| \leqslant 0.5] $ which is zero in the case of the law given by Figure~\ref{LawXExample}, then \begin{align*}
    P(C)&= 1/2, \\
    \mathrm{Cov}(X\mathbf{1}_C,X) &=1/16, \\
    \mathrm{Var}(X|C)&= 1/8, 
\end{align*} which imply that $(\sigma^{\nabla})^2 =11/48 \simeq 0.229 $ and $ e_P = 11/14 \simeq 0.786 $. Figure~\ref{Figdistribution_ztest_infoaux_H0} represents the distribution of $ \sqrt{n} \overline{X}_n,\sqrt{n} \overline{X}_n^\nabla $ which are respectively close to $ \mathcal{N}(0,\sigma^2) $ and $\mathcal{N}(0,(\sigma^\nabla)^2) $.
\includegraphics[width=260pt]{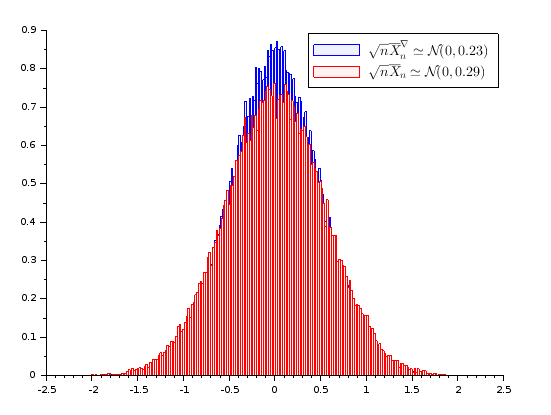}
  \captionof{figure}{Distribution of $ \sqrt{n}X_n $ and $ \sqrt{n}X_n^{\nabla } $}\label{Figdistribution_ztest_infoaux_H0}
  Figure~\ref{Figdistribution_ztest_infoaux_H1} displays the distribution of $ Z_n, Z_n^\nabla $ under $ (H_1) $ when $ \theta=0.5/\sqrt{n} $ and these statistics are asymptotically close to respectively $ \mathcal{N}(-\sqrt{6/7},1) $ and $ \mathcal{N}(-\sqrt{12/11},1) $.
\includegraphics[width=260pt]{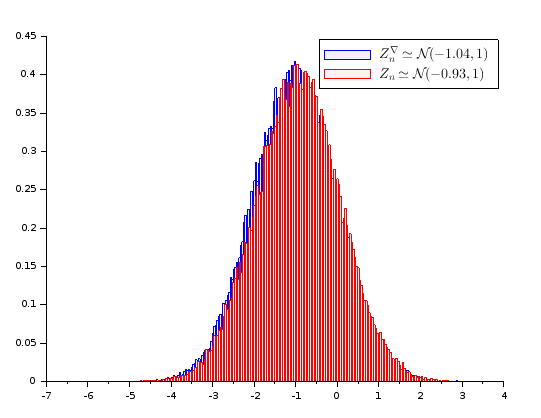}
  \captionof{figure}{Distribution of $ Z_n $ and $ Z_n^{\nabla } $ for $ n =100 $ under $ (H_1) $}\label{Figdistribution_ztest_infoaux_H1}

\subsubsection*{Conditional mean auxiliary information for the chi-square test} 
\textbf{General case.} Same notations of the previous paragraph are repeated. Suppose again that one exact auxiliary information is known but an uncertain auxiliary information can again also be considered.  In our case, $I=1 $, $ J_1=1 $ and \begin{align*}
    \Theta &= P[\mathcal{A}^*], \quad \widehat{\Theta} = \mathbb{P}_n[\mathcal{A}^*], \\
    \widetilde{\mathcal{B}} &= \mathcal{B}=\mathbb{E}[X|C], \quad \widehat{\mathcal{B}} = \mathbb{P}_n(X|C),
\end{align*} where $ \mathbb{P}_n(X|C) $ is the conditional expectation given by~\eqref{DefPnX|C}. The notation with $ \Theta $ is taken from the original paper and should not be confused with that of Section~\ref{Sec_ResultChiTest}. With these values, one have 
\begin{align*}
    K_{11} &= \mathrm{Var}(\widehat{\Theta}) = \mathrm{Var}(\mathbb{P}_n[\mathcal{A}^*]), \\
    K_{12} &= \mathrm{Cov}(\widehat{\mathcal{B}}, \widehat{\Theta}) = \mathrm{Cov}(\mathbb{P}_n(X|C), \mathbb{P}_n[\mathcal{A}^*]) \\
    &=\left( \mathrm{Cov}(\mathbb{P}_n(X|C), \mathbb{P}_n(A_i)) \right)_{1 \leqslant i \leqslant M-1}, \\
    K_{22}' &= \mathrm{Var}(\widehat{\mathcal{B}}) = \mathrm{Var}(\mathbb{P}_n(X|C)), \\
    K_{22}''&=\mathrm{Var}(\widetilde{\mathcal{B}})=0, \\
    K_{22} &= K_{22}'+K_{22}'' = K_{22}' = \mathrm{Var}(\mathbb{P}_n(X|C)).
\end{align*} Notice that values $ K_{22}',K_{22}'' $ and $ K_{22} $ does not change from the previous paragraph since the auxiliary information is exact and represents the same information. The elements $ K_{11}, K_{12}, K_{22} $ are still unknown or could not be expressed simply since the estimator $ \mathbb{P}_n(X|C) $ is a quotient of the empirical measure. Thus, the first suggested estimator of Tarima and Pavlov \begin{align*}
     \widehat{\Theta}^0 &= \mathbb{P}_n[\mathcal{A}^*]-K_{12}K_{22}^{-1} (\mathbb{P}_n(X|C)-\mathbb{E}[X|C]),
\end{align*} which exploits the auxiliary information $ \mathbb{E}[X|C] $ is impossible to evaluate. Values $ K_{11},K_{12},K_{22} $ can be estimated by the following consistent estimators $ \widehat{K}_{11} = \Sigma_{1,n} $ and \begin{align*} 
	\widehat{K}_{22} &= \frac{1}{n} \left(\mathbb{P}_n(X^2|C)-\mathbb{E}[X|C]^2 \right) \\
	\widehat{K}_{12} &= \frac{1}{n}  \left( \mathbb{P}_n(X \mathbf{1}_{A_i} | C)-\mathbb{E}[X|C] \mathbb{P}_n(A_i) \right)_{1 \leqslant i \leqslant M-1}.
\end{align*}With these estimators, Tarima and Pavlov suggest to use the statistic \begin{align*}
    \widehat{\Theta}^* &= \widehat{\Theta} - \widehat{K}_{12} \widehat{K}_{22}^{-1} (\widehat{\mathcal{B}}-\widetilde{\mathcal{B}}) \\
    &=\mathbb{P}_n[\mathcal{A}^*] - \widehat{K}_{12} \widehat{K}_{22}^{-1} (\mathbb{P}_n(X|C)-\mathbb{E}[X|C]).
\end{align*} By taking $ a_n= \sqrt{n} $ the first Tarima and Pavlov conditions mentioned in Section 1.3 of their paper are respected. More precisely, $ \zeta_{1n}=0, \Sigma_{22}''=0 $ and \begin{align*}
    \xi_n = a_n(\widehat{\Theta}-\Theta)&\underset{n\to+\infty}{\overset{\mathcal{L}}{\longrightarrow}} \xi = \mathcal{N}(\mathbf{0}_{M-1},\Sigma_{11}), \\
    \tau_n = a_n (\widehat{\mathcal{B}}-\mathcal{B}) &\underset{n\to+\infty}{\overset{\mathcal{L}}{\longrightarrow}} \tau = \mathcal{N}(0,\Sigma_{22}'), \\
    n K_{22}' = \mathrm{Var}(\alpha_n(X|C))  &\underset{n \to +\infty}{\longrightarrow} \Sigma_{22}',
    \end{align*} with $ \Sigma_{22}'= \mathrm{Var}(X|C)/P(C) $ where $ \mathrm{Var}(X|C) $ defined by~\eqref{DefVar|C} and $ \Sigma_{11}=\Sigma $ given by~\eqref{DefSig}. By Proposition 1 of Tarima and Pavlov, \begin{align*}
        a_n(\widehat{\Theta}^0-\Theta)
        \underset{n \to +\infty}{\overset{\mathcal{L}}{\longrightarrow}}& \mathcal{N}(\mathbf{0}_{M-1},\Sigma^\nabla),
    \end{align*} where \begin{align*}
        \Sigma^{\nabla} &= \Sigma_{11}-\Sigma_{12}\Sigma_{22}^{-1}\Sigma_{12}^t, \\
        \Sigma_{12} &= \mathrm{Cov}(\xi,\tau)= \frac{1}{P(C) } \left(\mathrm{Cov}(X\mathbf{1}_C,\mathbf{1}_{A_i})\right)_{1 \leqslant i \leqslant M-1}, \notag\\
        \Sigma_{22} &= \Sigma_{22}' = \frac{\mathrm{Var}(X|C)}{P(C)}. \notag
    \end{align*} Conditions of Proposition 2 of Tarima and Pavlov are respected since the following asymptotic behaviour is satisfied according to the distribution of large number: \begin{align*}
    a_n^2(\widehat{K}_{12}-K_{12}) = &a_n^2(\widehat{K}_{12}-\Sigma_{12}) - a_n^2(K_{12}-\Sigma_{12}) \\
    \underset{n\to+\infty}{\overset{\textrm{a.s.}}{\longrightarrow}} &\mathbf{0}_{M-1}, \\
    a_n^2(\widehat{K}_{22}-K_{22}) = &a_n^2(\widehat{K}_{22}-\Sigma_{22})-a_n^2(K_{22}-\Sigma_{22}) \\ \underset{n\to+\infty}{\overset{\textrm{a.s.}}{\longrightarrow}} &0.
\end{align*} By taking $ \mathbb{P}^{\nabla}_n[\mathcal {A}^*] =\widehat{\Theta}^* $, Proposition 2 of Tarima and Pavlov implies that \begin{align*}
    \sqrt{n}(\mathbb{P}^{\nabla}_n[\mathcal{A}^*]-P[\mathcal{A}^*])\underset{n \to +\infty}{\overset{\mathcal{L}}{\longrightarrow}}& \mathcal{N}(\mathbf{0}_{M-1},\Sigma^\nabla).
\end{align*} Matrix $ \Sigma^{\nabla} $ satisfies~\eqref{Cond1} then Proposition~\ref{PropChi2Stat} of this paper can be applied to the chi-square test which exploits the auxiliary information given by the knowledge of $ \mathbb{E}[X|C] $.\medskip
   
\textbf{Numerical simulation.} The previous results are applied with $ X $ distributed as Figure~\ref{LawXExample} and these following values: $ C = \{|X| \leqslant 0.5\}, M = 2, \mathcal{A}=\{A,A^C\} $ where $ A = \{X \leqslant 0\} $ satisfies $ P(A)=1/2 $.  With these values, the auxiliary information is given by $$ \mathbb{E}[X|C] = \mathbb{E}\left[X|\ |X|\leqslant 0.5\right] = 0 , $$that is the statistician knows the mean of the interest random variable when this last one is between -0.5 and 0.5. In this case, \begin{align*}
    \Sigma &=P(A)(1-P(A))=1/4,\\
    \Sigma_{12} &= \frac{\mathrm{Cov}(X\mathbf{1}_C,\mathbf{1}_A)}{P(C)} =-1/6, \\
    \Sigma_{22} &=\frac{\mathrm{Var}(X|C)}{P(C)}=1/4, \\
 \Sigma^{\nabla} &= \Sigma-\frac{\Sigma_{12}^2}{ \Sigma_{22}} = 5/36.
 \end{align*}  By Proposition~\ref{PropChi2Stat}, the Pitman's ARE $ e_P $ is $$ e_P = \Sigma^\nabla/\Sigma=  5/9. $$ Figure~\ref{Figdistribution_chi2test_infaux_H0} represents the distribution of $ \sqrt{n}(\mathbb{P}_n(A)-P(A)) $ and $ \sqrt{n}(\mathbb{P}_n^\nabla(A)-P(A)) $ for large value of $ n $. 
 
  \includegraphics[width=250pt]{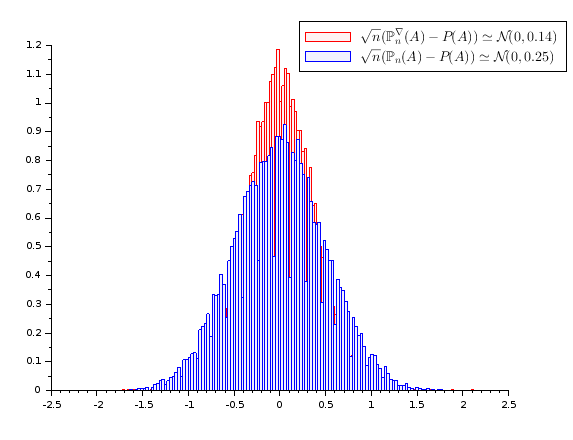}
  \captionof{figure}{Distribution of $ \sqrt{n}(\mathbb{P}_n(A)-P(A)) $ and $ \sqrt{n}(\mathbb{P}_n^\nabla(A)-P(A)) $}\label{Figdistribution_chi2test_infaux_H0}
  
 Figure~\ref{Figdistribution_chi2test_infaux_H1} represents the distribution function of $ \chi_n^2 $ and $\chi_n^{\nabla 2} $, for the hypothesis $ (H_1) $ with $ \mathbf{h}=0.5 $ and $ n=100$, which are respectively close to $ \chi^2(1;1) $ and $ \chi^2(1;9/5) $.

  \includegraphics[width=250pt]{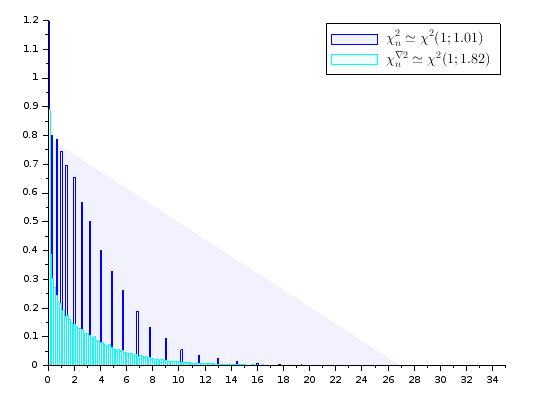}
  \captionof{figure}{Distribution of $ \chi_n^2 $ and $ \chi_n^{\nabla 2} $ for $ n = 100 $ under $ (H_1) $}\label{Figdistribution_chi2test_infaux_H1}

\section*{Acknowledgements}
I would like to sincerely thank the reviewers for their help which was invaluable to me. They gave me leads that I had not thought, for their advice which gave more consistency to this paper, made it easier to read and for the time they took to underline the mistakes.

\bibliographystyle{apalike}
\bibliography{biblio}

\end{multicols}

\end{document}